\documentclass[12pt,reqno]{amsart}
\usepackage{amsthm}
\usepackage{fullpage,url}
\usepackage[dvips]{hyperref}
 \usepackage{amscd}   % for commutative diagrams
\usepackage[all]{xy} % for complicated commutative diagrams
\usepackage{amssymb,amsmath}
\usepackage[dvips]{color}
\usepackage{setspace}
\DeclareFontEncoding{OT2}{}{} % to enable usage of cyrillic fonts

\newcommand{\EE}{\mathcal{E}}

\newcommand{\ZZ}{\mathbb{Z}}

\newcommand{\HH}{\mathbb{H}}
\newcommand{\FF}{\mathbb{F}}
\newcommand{\TT}{\mathbb{T}}

\newcommand{\RR}{\mathbb{R}}
\newcommand{\PP}{\mathcal{P}}
\newcommand{\pp}{\mathfrak{p}}
\newcommand{\LL}{\mathcal{L}}

\newcommand{\cH}{\mathcal{H}}
\newcommand{\br}{\bar{\rho}}

\newcommand{\nr}{\abs{\br}^2}

\DeclareMathOperator{\im}{Im}
\DeclareMathOperator{\re}{Re}
\DeclareMathOperator{\htt}{ht}
\newcommand{\e}[1]{\ar@{-}[#1]}
\newcommand{\ee}[1]{\ar@{=}[#1]}
\newcommand{\n}{*+[o][F-]{ }} 
\newcommand{\lul}[1]{\ar@{}[l]_<<{#1}}
\newcommand{\rul}[1]{\ar@{}[r]^<<{#1}}
\newcommand{\ldl}[1]{\ar@{}[l]^<<{#1}}
\newcommand{\rdl}[1]{\ar@{}[r]_<<{#1}}
\providecommand{\abs}[1]{\lvert#1\rvert}

\swapnumbers
\newtheorem{theorem}{Theorem}[section]
\newtheorem{lemma}[theorem]{Lemma}

\newtheorem{proposition}[theorem]{Proposition}

\theoremstyle{definition}

\theoremstyle{remark}

\newtheoremstyle{head}% name
{}% Space above
{}% Space below
{\bfseries}% Body font
{}% Indent amount (empty = no indent, \parindent = para indent)
{}% Thm head font
{}% Punctuation after thm head
{.5em}% Space after thm head: " " = normal interword space;
{}% Thm head spec (can be left empty, meaning `normal')

\theoremstyle{head}
\newtheorem{heading}[theorem]{}

\begin{document}

%\title[ql]{Reflection group of the quaternionic Lorentzian Leech lattice}
%
%\author{Tathagata Basak}
%\address{Department of Mathematics\\University of California, Berkeley\\Berkeley, CA 94720}
%\email{tathagat@math.berkeley.edu}
%\urladdr{http://www.math.berkeley.edu/\textasciitilde tathagat}
%
%\date{October 24, 2005}
%
%\begin{abstract}
% We find 14 reflections in the automorphism group of the the quaternionic Lorentzian Leech
% lattice $L$ that form the Coxeter diagram given by the incidence graph of projective plane
% over $\FF_{2}$. We prove that these 14 reflections generate the
% automorphism group of $L$. We find evidence that these reflections behave like
% the simple roots and the vector fixed by the diagram automorphisms behaves like
% the Weyl vector for the reflection group. 
%\end{abstract}
%
%\maketitle
%
{\bf {\large Reflection group of the quaternionic Lorentzian Leech lattice}}
\par
\vspace{.5cm}
Author : {\large Tathagata Basak}\par
\vspace{.5cm}
{\small Address : Department of Mathematics, University of California at Berkeley, Berkeley, CA 94720}
\par
{\small email : tathagat@math.berkeley.edu}\\

{\small Abstract: 
We study a second example of the phenomenon studied in the article ``The complex Lorentzian
Leech lattice and the bimonster''.
We find 14 roots in the automorphism group of the the quaternionic Lorentzian Leech
lattice $L$ that form the Coxeter diagram given by the incidence graph of projective plane
over $\FF_{2}$. We prove that the reflections in these 14 roots generate the
automorphism group of $L$. 
The investigation is guided by an analogy with the the theory of Weyl groups.
There is a unique point in the quaternionic hyperbolic
space fixed by the ``diagram automorphisms'' that we call the Weyl vector.
The unit multiples of the 14 roots forming the diagram are the analogs of the simple roots.
The 14 mirrors perpendicular to the simple roots are the mirrors that are closest to
the Weyl vector. \\
{\it keywords:} quaternionic Leech lattice, quaternionic hyperbolic reflection group, Y-groups, Weyl group, Coxeter diagram.

}
%
%******************************************************************************************************
%
\section{Introduction}
Let $\cH$ denote the ring of Hurwitz integers, consisting of the quaternions
$(a + b i + c j + d k)/2$ where $a,b,c$ and $d$ are integers, all congruent modulo 2.
Let $\Lambda$ and $E_8$ be the Leech lattice and $E_8$ root lattice respectively considered
as Hermitian lattices over $\cH$. Let $H$ be the 2 dimensional lattice $\cH \oplus \cH$ with
gram matrix $\bigl( \begin{smallmatrix} 0& 1+i \\ 1-i &0 \end{smallmatrix} \bigr)$
and let $L = \Lambda \oplus H$. Reflection groups of this and other quaternionic
Lorentzian lattices were studied by Allcock in \cite{dja:Leech} and \cite{dja:newcomplex}. 
\par 
In this article we investigate a second example of the phenomenon studied in the article \cite{tb:el}
on complex Lorentzian Leech lattice and bimonster. Here the base ring $\ZZ[e^{2\pi i/3}]$ is
replaced by the non-commutative ring $\cH$, the incidence graph of $P^2(\FF_3)$ is replaced by
$P^2(\FF_2)$ and there is a parallel story. For more on the bimonster and its presentations on
the $M_{666}$ and $\mbox{Inc}(P^2(\FF_3))$ diagrams, see \cite{cns:Y555}, \cite{cs:26} and
\cite{aai:geometryofsporadicgroups}. 
\par
We find 14 quaternionic reflections of order 4 in the reflection group of the Lorentzian
quaternionic Leech lattice $L$ that form the Coxeter diagram given by the incidence graph of 
the projective plane over $\FF_2$. (When I told D. Allcock about it, I got to know that he
had found this diagram too).
This 14 node diagram $D$ (see Fig. \ref{14}) is obtained by extending the $M_{444}$ diagram which 
comes up naturally, as the lattice $L$ is isomorphic to $3E_8 \oplus H$. (We shall describe an
explicit isomorphism from $\Lambda \oplus H $ to $3E_8\oplus H$ over $\cH$ later in \ref{calc},
because it is needed for our computations.)
The three hands of the $M_{444}$ diagram correspond to the three copies of $E_8$. 
\par
The main results are the following.
We see that there is a ``Weyl vector'' for $D$ in the quaternionic hyperbolic space that
is fixed by the diagram automorphisms and that the reflections in the 14 roots of $D$ generate the
reflection group of $L$. Allcock showed that the reflection group has finite index in the
automorphism group of $L$. We also see that the 14 mirrors of the reflections in the roots of $D$
are the mirrors closest to the Weyl vector, and in that sense, the roots of $D$ are the
analog of the simple roots.  
\par
The definitions and proof here are often similar to the ones in \cite{tb:el}, so we shall not
mention them at each step.
In view of the two examples, one would surely like to know whether there are other examples of similar
phenomenon for other Lorentzian lattices; and if a more conceptual meaning can
be given to these diagrams in a suitably general context that would explain all the numerical
coincidences seen in this and the previous article.
\par
\textit{Acknowledgments:} This paper was written as a part of my thesis
work at Berkeley. I would like to thank my advisor Prof. Richard Borcherds
for constant help, support and encouragement. I would like to thank Prof. Daniel Allcock for
generously sharing his ideas over many e-mails.
%
%******************************************************************************************************
%
\section{Preliminaries}
%
%******************************************************************************************************
%
\begin{heading}Notation\end{heading}
All the notations are borrowed from the article \cite{tb:el} on complex Lorentzian
leech lattice and are mostly consistent with the ones used in \cite{dja:Leech}.
As the ring of Eisenstein integers $\EE$ is replaced by the ring $\cH$ of Hurwitz
integers, the 14 node diagram $D$ replaces the 26 node diagram, $p$ replaces
$\theta = \sqrt{-3}$, $\FF_2$ replaces $\FF_3$, et cetera. 
\begin{tabbing}
$rad(D)$ X\= roots obtained by repeatedly reflecting of the roots in $D$
in each other.X\=\kill
$D$      \> The incidence graph of $p^2(\FF_2)$\\
$H$      \>  the hyperbolic cell over $\cH$ with Gram matrix 
$\bigl( \begin{smallmatrix} 
0&\bar{p} \\ p&0 
\end{smallmatrix} \bigr)$ \\
$\cH$     \> The ring of Hurwitz integers, a copy of the $D_4$ root lattice sitting inside the\\
\>quaternions.\\
$\HH$     \> The skew field of real quaternions.\\
$\Lambda$ \> The Leech lattice as a 6 dimensional negative definite $\cH$-lattice\\
$p$      \> $1-i$ \\
$\phi_r^{\alpha}$ \>  the $\alpha$-reflection in the vector $r$\\ 
$\xi$    \> $ (1 + i)/\sqrt{2}$\\
\end{tabbing}
%
%******************************************************************************************************
%
\begin{heading}The ring of quaternions\end{heading}
let $\cH$ be the ring of {\it Hurwitz integers} generated over $\ZZ$ by the 24 unit quaternions
$\pm 1, \pm i, \pm j, \pm k,$ and $(\pm 1 \pm i  \pm j \pm k)/2 $. The ring
$\cH $ consists of the elements $ (a + bi + cj + dk)/2 $, where $a, b, c, d$ are integers all congruent
modulo 2, with the standard multiplication rules $ i^2 = j^2 = k^2 = ijk = -1 $. When tensored with
$\mathbb{R}$ we get the skew field of quaternions called $\HH$. The conjugate
of $q = a+ bi + cj + dk$ is $\bar{q} = a - bi -cj -dk$. The {\it real part} of
$q$ is $\re (q) = a $ and the {\it imaginary part} is $\im (q) = q - \re(q)$. The {\it norm} of $q$ is 
$\abs{q}^2 = \bar{q}.q$.
\par
We let $\alpha = (1 + i + j + k)/2$ and $p = (1 - i)$.
In the constructions of the lattices given below, the number $p = 1 - i$
plays the role of $\sqrt{-3}$ in \cite{tb:el}.
$p$ generates a two sided ideal $\pp$ in $\cH$. We have
$i \equiv j \equiv k \equiv 1 \bmod \pp$. (Observe that $(1 - i)\alpha = 1 + j $,
so $j \equiv 1 \bmod \pp$). The group $\cH/\pp \cH$ is $\FF_4$, generated
by $0,1, \alpha, \bar{\alpha}$.
\par
The multiplicative group of units $\cH^*$ is $2\cdot A_4$. The quotient
$\cH^*/\lbrace \pm 1 \rbrace$ has four Sylow three subgroups generated by
$ \alpha, i\alpha, j\alpha, k\alpha$. The permutation representation of
$\cH^*/\lbrace \pm 1 \rbrace$ on the Sylow 3-subgroups identifies it with
the alternating group $A_4$.
\par
We also identify the quaternions $a + bi$ with the complex numbers. So any quaternion can be written as
$z_1 + z_2 j $ for complex numbers $z_1$ and $z_2$. The multiplication is defined by $j^2 = -1$
and $j z = \bar{z}j$. The complex conjugation becomes conjugation by the element $j$.
%
%******************************************************************************************************
%
\begin{heading} Lattices over Hurwitz integers \label{basis} \end{heading}
A general reference for lattices is \cite{cs:splag}.
An {\it $\cH$-lattice} is a free finitely generated right $\cH$-module with an $\cH$-valued bilinear form
$\langle\; ,\;\rangle$ satisfying $\overline{\langle x, y\rangle} = \langle y, x\rangle$,
$ \langle x, y\alpha\rangle = \langle x, y\rangle\alpha $ and
$\langle x\alpha, y\rangle =\bar{\alpha} \langle x, y\rangle$, for all $x, y$ in the lattice
and $\alpha$ in $\cH$. In this article, by a lattice we shall mean an $\cH$-lattice,
unless otherwise stated. Definite lattices will usually be negative definite.
The standard negative definite
lattice $\cH^n$ has the inner product $\langle x, y\rangle = -\bar{x}_1y_1 -\dotsb -\bar{x}_ny_n$.
The indefinite lattice given by $\cH \oplus \cH$ with inner product
$\langle (x,y), (x',y') \rangle = (\bar{x}, \bar{y})
\bigl( \begin{smallmatrix} 0 & \bar{p} \\ p & 0 \end{smallmatrix}\bigr)
\bigl( \begin{smallmatrix} x' \\ y' \end{smallmatrix} \bigr)$ is denoted by $H$.
We call $H$ the {\it hyperbolic cell}.
\par
Let $x$ be a vector in a $\cH$-lattice $K$. The norm $\abs{x}^2 = \langle x, x\rangle$ is
a rational integer. The $\ZZ$-module $K$ with the quadratic form given by the norm
will be called the {\it real form} of $K$. For example, the real form of the lattice
$\cH$ (with the norm multiplied by a factor of 2) is the $D_4$ root lattice.
\par
The $E_8$ root lattice can be  defined as a sub-lattice of $\cH^2$ as
$E_8 = \lbrace (x_1,x_2) \vert x_1 \equiv x_2 \bmod \pp \rbrace$. It has minimal norm $-2$ and the
underlying $\ZZ$-lattice is the usual $E_8$ root lattice. 
\par
The Leech lattice $\Lambda$ can be defined as a 6 dimensional negative definite $\cH$-lattice
with minimal norm $-4$ whose real form is the usual real Leech lattice. The automorphism group of
this lattice was studied by Wilson in \cite{raw:g24}. We quote the facts we need from there.
Let $\omega = (-1+i+j+k)/2$.
The lattice $\Lambda$ consists of all vectors  ($v_{\infty}, v_0, v_1, v_2, v_3,v_4)$ in $\cH^6$ such that
$v_2 \equiv v_3 \equiv v_4 \bmod \pp$,
$(v_1 + v_4)\bar{\omega} + (v_2 + v_3)\omega \equiv (v_0 + v_1)\omega + (v_2 + v_4)\bar{\omega} \equiv 0 \bmod 2$, 
and $-v_{\infty}(i+j+k) + v_0 + v_1 + v_2 + v_3 + v_4 \equiv 0 \bmod 2+2i$. The
inner product we use is $-1/2$ of the one used in \cite{raw:g24}, so that the following basis vectors have norm $-4$.
We use the following $\cH$-basis for the lattice $\Lambda$ given in \cite{raw:g24} for some computations.
\begin{align*}
bb[1] &= [2+2i , 0  ,  0  ,   0  ,  0  ,  0  ]\\    
bb[2] &= [  2  , 2  ,  0  ,   0  ,  0  ,  0  ]\\    
bb[3] &= [  0  , 2  ,  2  ,   0  ,  0  ,  0  ]\\    
bb[4] &= [i+j+k, 1  ,  1  ,   1  ,  1  ,  1  ]\\    
bb[5] &= [  0  , 0  ,1+k  , 1+j  ,1+j  ,1+k  ]\\    
bb[6] &= [  0  ,1+j ,1+j  , 1+k  , 0   ,1+k  ]\\    
\end{align*}
Let $L$ be the Lorentzian lattice $L = \Lambda \oplus H \cong 3E_8 \oplus H $. The real form of this
lattice is $II_{4,28}$.
$E_8$, $\Lambda$, $H$ and $L$ each satisfy $L'p = L$, where $L'$ is the \textit{dual lattice} of $L$
defined by $L' = \lbrace x \in L \otimes \HH : \langle x, y\rangle \in \cH \forall y \in L\rbrace$.
\par
%
%******************************************************************************************************
%
\begin{heading} Quaternionic reflections\end{heading}
A \textit{$\mu$-reflection} in a vector $r$ of a lattice is given by
\begin{equation}
\phi_r^{\mu}(v) = v - r (1 - \mu) \langle r,v\rangle/\abs{r}^2
\end{equation}
where $\mu \neq 1$ is a unit in $\cH$.
Note that 
\begin{equation}
\phi_{r\alpha}^{\mu} = \phi_{r}^{\alpha \mu \alpha^{-1}}
\end{equation}
for any unit $\alpha \in \cH$.
Note also that $\phi_r^{\mu}$ can be characterized as the automorphism of the lattice that fixes the
orthogonal complement of $r$ and multiplies $r$ by the root of unity $\mu$.
It follows that for every automorphism $\gamma$ of the lattice we have
$\gamma \phi_r^{\mu} \gamma^{-1} = \phi_{\gamma r}^{\mu}$.
We say the two reflections $\phi_1$ and $\phi_2$ \textit{braid} if $\phi_1 \phi_2 \phi_1 =  \phi_2 \phi_1 \phi_2$.
A \textit{root} of a negative definite or a Lorentzian lattice is
a lattice vector of negative norm such that there is a nontrivial reflection in it that is an
automorphism of the lattice. If $L'p  = L$, the roots of $L$ are all
the vectors of norm $-2$. For such a root $r$, $R(L)$ contains six reflections
of order four given by $\phi_r^{\pm i},\phi_r^{\pm j}$ and $\phi_r^{\pm k}$.
The square of each of them is the order 2 reflection $\phi_r^{-}$.
Note that the six units $\pm i, \pm j,$ and $\pm k$ form a conjugacy class in $\cH^*$.  
For a set of reflections (or roots) that either braid or commute, we form the \textit{Coxeter diagram}
by taking one vertex for each reflection and joining them only if the reflections braid.
\par
For computational purpose we note the following (one needs to be
careful because $\cH$ is not commutative). Let $e_1, \dotsb, e_n$ be a basis for a lattice.
For a linear transformation $\phi$, if $\phi(e_t) = \sum_s e_s \phi_{s t}$ then
let $mat(\phi) = ((\phi_{s t})) $; if a vector $x = \sum e_t x_t$ is written in
coordinates as a column vector then $\phi(x) = mat(\phi)(x_1,\dotsb,x_n)'$.
%
%******************************************************************************************************
%
\begin{heading} 
Computation for explicit isomorphism between $3E_8 \oplus H$ and $\Lambda \oplus H$ over $\cH$
\label{calc}
\end{heading}
The reflection group of the $\cH$-lattice $E_8$ can be generated by two
reflections of order 4 that braid with each other. So an $E_8$-diagram
for us, looks like an $A_2$ Dynkin diagram. We say that two diagrams
are orthogonal if the roots in the first diagram are orthogonal to the
roots in the second diagram. We need to find three orthogonal
$E_8$-diagrams in the lattice $\Lambda \oplus H$
and a hyperbolic cell orthogonal to this $3E_8$, to get an explicit change of basis matrix. 
We find the 6 vectors of the form $r = (l; 1 ,\bar{p}^{-1}(1 + \beta))$
( with $\beta \in \im{\HH}$ and $ l $ in the first shell of the leech lattice )
in $\Lambda \oplus H$ forming three copies of $E_8$
by a computer search. The steps of the computation are given below.
\par
First, note that $i$-reflections in two roots $r =(l; 1 ,\bar{p}^{-1}(1 + \beta))$ and 
$r' =(l'; 1 ,\bar{p}^{-1}(1 + \beta'))$ will commute if $\abs{l - l'}^2 = -4 $ and 
$\beta - \beta' = [l, l'] $, where $[l,l'] = \im \langle l, l' \rangle $.
The $i$-reflections in $r$ and $r'$ will braid if $\abs{l - l'}^2 = 6$ and 
$\beta - \beta' = [l,l'] \pm i $. 
\par
It is easy to find one copy of $E_8$. One just have to find two vectors $l_1, l_2$ in first shell of Leech
lattice that are at a distance $-6$ and find $\beta_1, \beta_2$ accordingly.
\par
Now we generate a large list of vectors in the first shell of Leech lattice using the basis $bb$ given
in \ref{basis}. We find all the vectors $l$ in the first shell (actually from the almost complete list that
we had), that might give a root $(l; 1,*)$ of an orthogonal $E_8$. 
The conditions that $l$ has to satisfy are $\abs{ l - l_1}^2 = \abs{l - l_2}^2 = -4 $ and
$[l_1,l_2] - [l_1,l] + [l_2, l] = i $. This gives a small list of vectors $l$.
(It is amusing to note that the equation looks like a ``co-cycle condition'').
\par
Next, find pairs of vectors $l_3, l_4$ from the previous list that can actually give an
$E_8$ orthogonal to the first one. The conditions to be satisfied are 
$[l_1,l_2]- [l_3,l_4]+ [l_2,l_4]- [l_1,l_3] = 0$ and $\abs{l_3 - l_4}^2 = - 6$.
This way one can find two $E_8$ diagrams orthogonal to the first one, that happen to
be orthogonal to each other too. Thus we get the 6 vectors forming a $3E_8$ diagram.
\par
Now we take the orthogonal complement and find two norm zero vectors in the complement forming
an hyperbolic cell.  
\par
The rows of the matrix given below are the actual coordinates for the 8 vectors found
by above calculation. The first six root vectors form a basis for the lattice $3 E_8$ and the last
two norm zero vectors form a basis of a hyperbolic cell orthogonal to the first six vectors.
\begin{scriptsize}
\begin{equation*}
\begin{pmatrix}
%1
2 & 2 & 0 & 0 & 0 & 0 & 1 & 1 \\
%2
i+j+k & 1 & 1 & 1 & 1 & 1 & 1 & 2 + i + j \\
%3
\tfrac{1+i+j+k}{2} & \tfrac{3+i-j+k}{2}& \tfrac{1-i+j-k}{2} & \tfrac{1-i-j+k}{2} 
&\tfrac{-1-i-j-k}{2}& \tfrac{1-i-j+k}{2}& 1                & \tfrac{3+i+j+k}{2} \\
%4
1+i+k & 1 & k & -k & 1 & -j & 1 & \tfrac{3+i+j+k}{2} \\
%5
1+i+k & 1 & i  & -j & -i & 1 & 1 & \tfrac{3+i+j+k}{2} \\
%6
1 + i & 1 + k & 1 - i & 1  -k & 0 & 0 & 1 & \tfrac{3+i+j+k}{2} \\
%7
\tfrac{-7+ i -3 j -5 k}{2} & \tfrac{-7+3 i -3 j +k}{2} & \tfrac{-1 + 3 i - j + k}{2} & \tfrac{-1+i - j + 3 k}{2}
&\tfrac{-1 + 3 i + j + 3 k}{2} & \tfrac{-3 + i -j + k}{2}  & \tfrac{-5 + 3 i - 3 j + k}{2}& -4  +2 i -3 j  - k  \\
%8
\tfrac{-5 - 3 i -j - 5k}{2} & \tfrac{-7 + i -j -3 k}{2}, &  \tfrac{-1+i-j-k}{2} & \tfrac{-3 + i-j+k}{2} 
&\tfrac{-1 + 3 i-j+k}{2} & \tfrac{-3 + i+j-k}{2} & \tfrac{-5 + i-j-k}{2} & -4  -2 j - 2 k
\end{pmatrix}
\end{equation*}
\end{scriptsize}
%
%the same vectors copy pasted from the actual program in raw form
%\begin{tiny}
%\begin{verbatim}
%1 [2, 0; 2, 0; 0, 0; 0, 0; 0, 0; 0, 0; 1, 0; 1, 0]
%2 [I, 1 + I; 1, 0; 1, 0; 1, 0; 1, 0; 1, 0; 1, 0; 2 + I, 1]
%3 [1/2+1/2*I,1/2+1/2*I ; 3/2+1/2*I,-1/2+1/2*I ; 1/2-1/2*I,1/2-1/2*I ; 1/2-1/2*I,-1/2+1/2*I;
%  -1/2-1/2*I,-1/2-1/2*I; 1/2-1/2*I,-1/2+1/2*I ; 1, 0                ; 3/2+1/2*I, 1/2+1/2*I]
%4 [1 + I, I; 1, 0; 0, I; 0, -I; 1, 0; 0, -1; 1, 0; 3/2 + 1/2*I, 1/2 + 1/2*I]
%5 [1 + I, I; 1, 0; I, 0; 0, -1; -I, 0; 1, 0; 1, 0; 3/2 + 1/2*I, 1/2 + 1/2*I]
%6 [1 + I, 0; 1, I; 1 - I, 0; 1, -I; 0, 0; 0, 0; 1, 0; 3/2 + 1/2*I, 1/2 + 1/2*I]
%7 [-7/2+1/2*I,-3/2-5/2*I ; -7/2+3/2*I,-3/2+1/2*I ; -1/2+3/2*I,-1/2+1/2*I ; -1/2+1/2*I,-1/2+3/2*I;
%   -1/2+3/2*I,1/2 +3/2*I ; -3/2+1/2*I,-1/2+1/2*I ; -5/2+3/2*I,-3/2+1/2*I ; -4  +2*I  ,-3  -I]
%8 [-5/2-3/2*I,-1/2-5/2*I ; -7/2+1/2*I,-1/2-3/2*I ; -1/2+1/2*I,-1/2-1/2*I ; -3/2+1/2*I,-1/2+1/2*I;
%   -1/2+3/2*I,-1/2+1/2*I ; -3/2+1/2*I, 1/2-1/2*I ; -5/2+1/2*I,-1/2-1/2*I ; -4        , -2 - 2*I]
%\end{verbatim}
%\end{tiny}
%
%*************************************************************************************************
%
\section{the 14 node diagram}
\begin{figure}
\centerline{
\xymatrix@1@=16pt@!{
&\n \ldl{d_1}             &  &  &  &  &  &                   \n \rdl{d_2}& \\
& & \n \ldl{c_1}   \e{ul}     &  &  &  &     \n \rdl{c_2}   \e{ur}  &    &  \\
& & & \n \ldl{b_2}   \e{ul}      &  &        \n \rdl{b_2} \e{ur}  & &    &  \\
& & & & \n \rdl{a} \e{ul} \e{ur} \e{d}                          & & &    &  \\ 
& & & & \n \rdl{b_3}    \e{d}                                   & & &    &  \\ 
& & & & \n \rdl{c_3}    \e{d}                                   & & &    &  \\
& & & & \n \rdl{d_3}                                            & & &    &
}
\xymatrix@1@=16pt@!{
&             &             &\n \e{rrrr} \lul{d_3} &             &             &             &\n \e{ddd} \rul{e_2}   &\\
&             &             &\n \e{u} \e{rrddd}\ldl{c_3}    &        &             &\n \e{ur} \e{llllddd}\rdl{b_2}&             &\\
&\n\e{uurr} \lul{e_1}&      &             &\n \e{lll}\e{d}\rul{b_1} &\n \e{l}\e{dd}\e{ur}\rdl{a}&         &             &\\
&             &    &\n\e{uu}\e{dl}\e{r}\lul{f}&\n \e{rrr}\rdl{c_1}&     &             &\n \e{ddll}  \rdl{d_1}        &\\
&             &\n \e{dl}\e{ur}\rdl{c_2}&           &             &\n \e{d} \rul{b_3}&             &             &\\
&\n\e{uuu} \ldl{d_2} &       &            &             &\n \e{llll} \rdl{e_3}&             &           &
}
}
\caption{The diagrams $M_{444}$ and $D$ }
\label{14}
\end{figure}
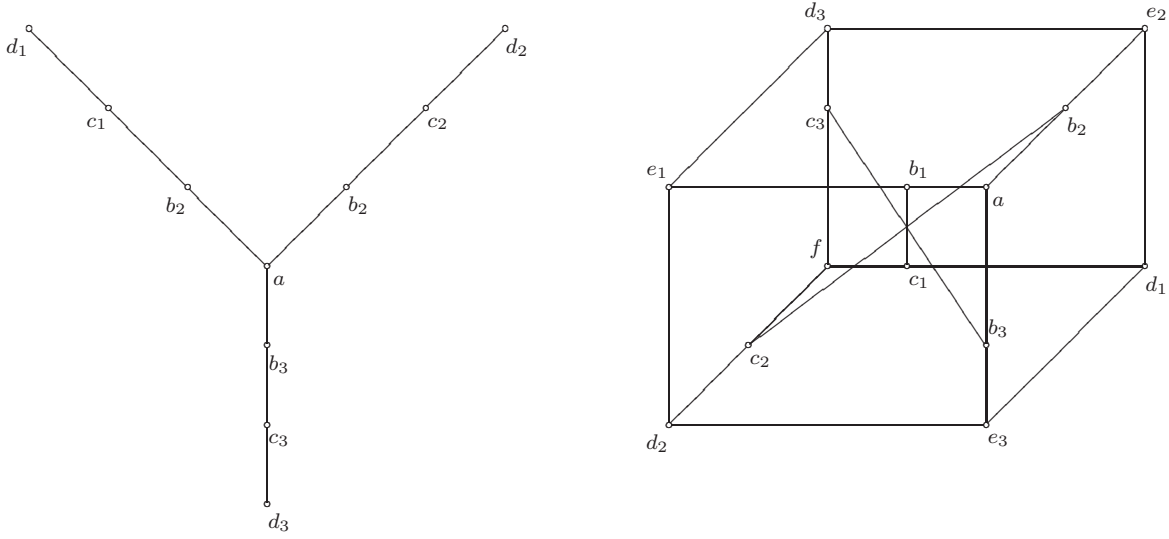

%
%******************************************************************************************************
%
\begin{heading} the diagram of 14 roots\end{heading}
The $i $-reflection of order 4 in  a root $r$ braids with the $i$-reflection
in a root $r'$ if $\langle r, r'\rangle = p $.  
In this section we work in the co-ordinates $3E_8 \oplus H$ for $L$.
We can find 10 roots $a, b_s, c_s, d_s$ for $ s = 1, 2,3 $ forming an $M_{444}$ diagram
inside the reflection group $ R(L)$.
See \cite{aai:geometryofsporadicgroups} for more on these groups (called $Y$-groups there).
The  hands of the $M_{444}$-diagram correspond to the three copies of $E_8$ in $L$, in the sense
that $c_s, d_s$ generate the reflection group of the $s$-th $E_8$, the $b_s$ are the affinizing node
and $a$ is the hyperbolizing node.
These 10 roots can be extended to a set of $14$ roots forming the incidence graph of
$\mathbb{P}^2(\FF_2)$. The roots $a,c_s, e_s$ correspond to the points of $\mathbb{P}^2(\FF_2)$
(See Fig. \ref{14})
and $f, b_s, d_s$ correspond to the lines. There is an edge between them if the point lies on the line.
This 14 vertex diagram is called $D$. the seven roots $a,c_s, e_s$ (or their unit multiples)
are called ``points'' and the seven roots $f, b_s, d_s$ (or their unit multiples) are called ``lines''.
%******************************************************************************
%******************************************************************************
%******************************************************************************
%******************************************************************************
One choice of the explicit coordinates of these roots are given by
\begin{equation*}
\begin{matrix}
r[1  ]  = & a   =[&          &,&               &;&   1 &,&  -1 &]\\ 
r[1+s]  = & c_s =[&(1,-1)_s  &,&               &;&     &,&     &]\\ 
r[4+s]  = & e_s =[&          &,& (1,1)_{t u}    &;&  -i &,&  -1 &]\\ 
r[8]    = & f   =[& (\;\;,p)^3 & &               &;& -ip &,&  -p &]\\ 
r[8+s ] = & b_s =[&(\;\;,p)_s  &,&               &;&  -1 &,&     &]\\ 
r[11+s] = & d_s =[&(-p,\;)_s &,&               &;&     &,&     &]
\end{matrix}
\end{equation*}
An $E_8$ vector with subscript $s$ means we put it in place of the $s$-th $E_8$, while
the subscript $t u$ means that we put it at the $t$-th and $u$-th place. The indices $s, t$ and $u$
are in cyclic permutation of $(1,2,3)$. Blank spaces are to be filled with zero. 
%
%******************************************************************************************************
%
\begin{heading}Linear relations among the roots of $D$  \label{rel} \end{heading}
Let $x$ be a point and $l$ be a line of $D$. $\langle x,l \rangle$ is equal to $p$ or $0$ depending
on whether there is a edge between $x$ and $l$ or not. Using this we find that the vector 
\begin{equation} 
w_{\PP}  = l\bar{p} + \sum_{x \in l}x
\label{wp}
\end{equation}
is perpendicular to the points
and has norm $2$. But there is only one such vector. Hence we get the relations 
$l\bar{p} + \sum_{x \in l}x = l'\bar{p} + \sum_{x \in l'}x$.
for any two lines $l$ and $l'$. These relations generate all linear relations in the 14 roots.
Using the automorphism $\sigma$ that takes $l$ to $x$ and $x$ to $li$ we see that the element
\begin{equation} 
w_{\LL} = x p + \sum_{x \in l} l
\label{wl}
\end{equation}
is constant for all point $x$. We shall see that
the elements $w_{\PP}$ and $w_{\LL}$ determine points in the quaternionic hyperbolic space fixed by
the diagram automorphisms $PGL_2(\FF_2)$.
%
%*************************************************************************************************
%
\section{Generators for the reflection group of L}
%
%******************************************************************************************************
%
\begin{heading}The Heisenberg group\end{heading}
We follow the definitions and notations of section 6 of the article \cite{dja:Leech}. 
Let $\TT$ be the \textit{Heisenberg group} generated by the translations $T_{\lambda,z}$ for every $\lambda $ in
$\Lambda = \Lambda \cap \Lambda'p$. In section 6
Allcock mentions that formulae $(3.2)$ - $(3.5)$ of \cite{dja:Leech} holds in the quaternionic case too.
We need to do a little more calculation to get a little stronger version of Theorem 6.1 of \cite{dja:Leech}.
We use the three roots $r_1 = (0^6; 1, i)$, $r_2 = (0^6; 1, -1)$ and $r_3 = ( 0^6; 1, -\epsilon)$ where
$\epsilon = ( 1 - i - j + k)/2 $. 
%
%******************************************************************************************************
%
\begin{lemma} Let $\epsilon = ( 1- i - j + k)/2$, $r_3$ and $r_2$ be the roots $(0^6;1,-\epsilon)$
and $(0^6;1, -1)$ and let $R = \phi_{r_3}^i \phi_{r_2}^j $. Then 
\begin{equation}
T_{\lambda,z}^{-1} R T_{\lambda,z} R^{-1} = T_{\lambda(\bar{\epsilon} - 1), \epsilon z \bar{\epsilon} - z +
\im\langle \lambda \bar{\epsilon} , -\lambda\rangle }
\label{commutator}
\end{equation}
Since $(\bar{\epsilon} - 1)$ is an unit of $\cH$ we see that for every $\lambda \in \Lambda$, $R(L)$ contains
a reflection in $\lambda$.
\label{conjwitht}
\end{lemma}
\begin{proof}
$R$ is of the form $\bigl( \begin{smallmatrix} I & 0 \\ 0 & R_H \end{smallmatrix}\bigr)$ where
$R_H = \bigl( \begin{smallmatrix} \epsilon & 0 \\ u  & \delta \end{smallmatrix}\bigr)$ is the matrix acting
on $H$ with $u = (3 - i + j - k)/2$ and $\delta \bar{p}^{-1} = \bar{p}^{-1} \epsilon$.
Now $R^{-1}_H = \bigl( \begin{smallmatrix} \epsilon^{-1} & 0 \\ v  & \delta^{-1} \end{smallmatrix}\bigr)$
where $v = -\delta^{-1}u \epsilon^{-1}$.
\par
Now just multiply the matrices to see that 
\begin{equation*}
\begin{split}
R T_{\lambda, z} R^{-1} 
&= \begin{pmatrix} I &         & \\ & \epsilon & \\  & u & \delta \end{pmatrix}
  \begin{pmatrix} I &\lambda  & \\ & 1        & \\ -{\bar{p}}^{-1}\lambda^*  &\bar{p}^{-1}(z - \lambda^2/2) & 1 \end{pmatrix}
   \begin{pmatrix} I &         & \\ &\bar{ \epsilon} & \\  & v &\bar{ \delta} \end{pmatrix}\\
&=\begin{pmatrix} I &\lambda  & \\ & \epsilon & \\-\delta{\bar{p}}^{-1}\lambda^* &u +\delta \bar{p}^{-1}(z - \lambda^2/2) & \delta \end{pmatrix}
  \begin{pmatrix} I &         & \\ & \bar{\epsilon} & \\  & v & \bar{\delta} \end{pmatrix}\\
&=\begin{pmatrix} I &\lambda\epsilon^{-1}  & \\ & 1 & \\-\delta{\bar{p}}^{-1}\lambda^*  & x & 1 \end{pmatrix}
\end{split}
\end{equation*}
where 
\begin{equation*}
x = (u + \delta{\bar{p}}^{-1}(z - \lambda^2/2))\bar{\epsilon} + \delta v
  = \delta{\bar{p}}^{-1}(z - \lambda^2/2)\bar{\epsilon} 
  ={\bar{p}}^{-1}\epsilon(z - \lambda^2/2)\bar{\epsilon}  
  ={\bar{p}}^{-1}(\epsilon z \bar{\epsilon} - \lambda^2/2)
\end{equation*}
In the second equality we use $u\bar{\epsilon}  = -\delta v$ and in the third $\delta {\bar{p}}^{-1} = {\bar{p}}^{-1}\epsilon$.
Also note that $-\delta \bar{p}^{-1} \lambda^* = -\bar{p}^{-1}\epsilon\lambda^* = -\bar{p}^{-1}(\lambda\bar{\epsilon})^*$.
Thus we have
\begin{equation}
R T_{\lambda, z} R^{-1} = T_{\lambda \bar{\epsilon}, \epsilon z \bar{\epsilon }}
\label{RTRinv}
\end{equation}
Now \eqref{commutator} follows from \eqref{RTRinv} and
$T_{\lambda,z} T_{\lambda',z'} = T_{\lambda + \lambda', z + z' + \im\langle \lambda',\lambda\rangle}$ (equation (3.2) in \cite{dja:Leech}).
\end{proof}
%
%******************************************************************************************************
%
\begin{lemma}
Let $\Psi$ be the set of roots that are unit multiples of the roots of the form $(\lambda; 1 , *)$.
The reflections in $\Psi$ act transitively on all roots. The reflections in all the roots of the
form $(\lambda, 1, *)$  generate the reflection group of $L$.
\end{lemma}
\begin{proof}
The calculation here is almost identical to Theorem 6.2 in \cite{dja:Leech} which uses the idea in \cite{jhc:automorphism}.
For this lemma only, let $h(\lambda; \mu, \eta) =  \abs{\mu}$. (Later we are going to use a different definition of
height). The roots $r$ in $\Psi$ are the ones with $h(r) = 1$. 
We show that if we have a root $r$ with $h(r) > 1$, then we can $i$-reflect it in a root of $\Psi$
to decrease its height. The covering radius of the Leech lattice is used in this proof and it has just the
right value to make things work. 
\par
Let $y = (l; 1, \bar{p}^{-1}( \alpha - l^2/2))$ be a multiple of a root $r$ with $h(r) > 1$.
We have $\abs{y}^2 \in (-2,0) $, which amount to $ \re (\alpha) \in (-1,0)$ because $\abs{y}^2 = 2 \re (\alpha)$.
Consider the $i$-reflection in the root $r = (\lambda; 1, \bar{p}^{-1}(-1 -\lambda^2/2 + \beta + n)$; where
$\beta \in \im\HH$ is determined so that $\bar{p}^{-1}(-1 -\lambda^2/2 + \beta) $ is in $\cH$ and $n$
( to be chosen later)
can be any element of $\im\pp = \pp \cap \im\HH$. Calculation yields $\langle r, y \rangle = -2(a + b) $, where 
\begin{equation*}
-2a = -\frac{1}{2}\abs{l - \lambda}^2 - 1 + Re(\alpha) \in \RR 
\end{equation*}
and 
\begin{equation*}
-2b = \im(\alpha) + \im\langle \lambda, l\rangle - \beta - n \in \im\HH
\end{equation*}
So $h(\phi_r^i(y)) = \bar{p}( 1 - (1 - i)(a + b))$. Thus we want to make
$\abs{ 1 - (1 - i)(a + b)}^2 < 1 $, which amount to 
\begin{equation}
\abs{1/2 - a}^2 + \abs{i/2 - b}^2 < 1/2 
\label{lthalf}
\end{equation}
Because the covering norm of Leech lattice  is 2 we can make $\abs{l - \lambda}^2 \in [-2,0] $.
This, together with $Re(\alpha) \in ( - 1, 0) $ gives $ a \in (0,1)$. So $\abs{1/2 - a}^2 <1/4$.
As for the second term of \eqref{lthalf}, $(i -2b)$ is in $\im\HH$, and in the expression for $-2b$
we are free to choose $n \in \im{\pp}$ which forms a copy of $D_3$ root lattice:
$\lbrace(a i + bj+ ck) : a + b + c \equiv 0 \bmod 2\rbrace$ in $\im\HH$. The covering radius of
$D_3$ is 1. So we can make the norm of $(i -2b)$ less than 1 by choice of $n$ and thus
make $\abs{i/2 - b}^2 \leq 1/4$. 
\par
So if $\phi_r^\mu$ is a reflection in any root $r$, after conjugating finitely many times
by $i$-reflections in roots of the form $ ( \lambda; 1, *)$ we get a reflection $\phi_{r'}^\zeta$ in a root
$r'= r'' u$ where $r''= (\lambda''; 1, *)$ and  $u$ is an unit. But then $\phi_{r'}^{\zeta} = \phi_{r''}^{u \zeta u^{-1}} $. 
Thus $\phi_r^{\mu}$ can be obtained as a product of reflections in the roots of the form $(\lambda; 1, *)$.
\end{proof}
%
%******************************************************************************************************
%
\begin{lemma}
Let $\lambda_1, \lambda_2, \dotsb, \lambda_{24}$ be elements of $\Lambda$ that make a $\ZZ$-basis.
Let $r_1 = (0^6; 1, i) $ and $r_2$, $r_3$ be as given in Lemma \ref{conjwitht}.
Fix $z_s$ such that $T_{\lambda_s, z_s} \in \TT$, for $s = 1, \dotsb, 24$.
Let $R_1$ temporarily denote the group generated by all the reflections in the 81 roots 
$T_{\lambda_s, z_s}(r_t)$, $T_{0,i+j} (r_t)$, $T_{0,i+k} (r_t)$ 
and $r_t$ where $s = 1, \dotsb, 24$ and  $ t = 1, 2, 3$.
Then $R_1$ contains the Heisenberg group $\TT$. In fact $R_1$ is equal to the reflection group of $L$.
\label{gen}
\end{lemma}
\begin{proof}
From \eqref{commutator} in lemma \ref{conjwitht} we get that for each $\lambda_s$, the group $R_1$ contains a
translation in the vector $\lambda_s(\bar{\epsilon} - 1)$. These vectors form a $\ZZ$-basis of $L$ as
$\bar{\epsilon} - 1$ is an unit. Using
$T_{\lambda,z}\circ T_{\lambda',z'} = T_{\lambda + \lambda', z + z' + \im\langle \lambda',\lambda\rangle }$
and $T_{\lambda,z}^{-1} = T_{-\lambda, -z}$ (equations (3.2) and (3.3) in \cite{dja:Leech})
we see that $R_1$ contains translation in every vector of $\Lambda$.
\par
Now we argue that all the central translations of the form $T_{0,z}$ are in $R_1$.
Choosing $\lambda$ and $\lambda'$ such that $\langle \lambda', \lambda\rangle = p$ and using the identity 
$T_{\lambda,z}^{-1}T_{\lambda',z'}^{-1}T_{\lambda,z}T_{\lambda',z'} = T_{0,2\im\langle \lambda',\lambda\rangle}$
(equation (3.4) in \cite{dja:Leech}) conclude that the central translation $T_{0,2i}$ is in $R_1$.
Similarly taking $\langle \lambda', \lambda\rangle $ to equal $p\alpha = (1 + j)$ and $p( 1 + i - j + k)/2 = (1 + k)$
respectively, it follows that the central translation $T_{0,2j}$ and $T_{0,2k}$ are also in $R_1$.
From \eqref{RTRinv} in \ref{conjwitht}, it follows that 
$T_{0,z}^{-1} R T_{0,z} R^{-1} = T_{0, \epsilon z \bar{\epsilon} - z} $ which equals 
$T_{0,-i-2j-k} $ for $z = i+j$ and $T_{0,-j-k}$ for $z = i+k$. So these central translations are in $R_1$ too.
We found that $T_{0,i+k}$ and $T_{0,j+k}$ are in $R_1$, as are $T_{0,2i}, T_{0,2j}$ and $ T_{0,2k}$.
These central translations clearly generate all the central translations of the form $T_{0,z}$, $z\in \pp\cH $.
So $R_1$ contains $\TT$.
\par
the orbit of $(0^6;1, -1)$ under $\TT$ is all roots of the form
$ ( \lambda, 1, \bar{p}^{-1}( \beta - 1 - \lambda^2/2))$.
So all these roots are in $R_1$ and we have already seen that these generate
the whole reflection group of $L$. 
\par
For the actual computation multiply the basis vectors $bb$ given in \ref{basis}
by $1, i, j, (-1+i+j+k)/2$ to get a basis of $\Lambda$ over $\ZZ$.
\end{proof}
%
%*************************************************************************************
%
\section{The fixed points of diagram automorphisms and height of a root}
%
%******************************************************************************************************
%
\begin{heading} The fixed points under diagram automorphism\end{heading}
The group $PGL_2(\FF_2)$ acts on the diagram $D$ and this induces a linear action of $PGL_2(\FF_2)$ on $L$.
The graph automorphism switching points with lines in $Inc(\mathbb{P}^2(\FF_2))$ lifts to a
automorphism $\sigma$ taking a line $l$ to a point $x$ and  $x$ to  $li $
(note: $\langle x, l \rangle = p$ if and only if $ \langle li, x \rangle= p$).
This gives action of the extended diagram automorphism group
$Q = 8\cdot PGL_2(\FF_2)$ on $L$ and hence on the \textit{quaternionic hyperbolic space} $\HH H^7$,
which consists of the positive norm lines in $L \otimes \RR$. The vectors
$w_{\PP}$ and $w_{\LL}$ defined in \eqref{wp} and \eqref{wl} of section \ref{rel}
span the 2-dimensional space fixed by the action of $PGL_2(\FF_2)$.
From \eqref{wp} and \eqref{wl} it follows that $\sigma(w_{\pp}) = w_{\LL}i$ and $\sigma(w_{\LL}) = w_{\PP}$.
\par
Let $\Sigma_{\PP}$ and $\Sigma_{\LL} $ be the sum of the points and the sum of the lines respectively.
These too are fixed by the $PGL_2(\FF_2)$ action and $\sigma$ takes $\Sigma_{\LL}$ to $\Sigma_{\PP}$ and
$\Sigma_{\PP}$ to $\Sigma_{\LL}i $. 
So, there is a unique fixed point in $\HH H^7$ under the action of this extended group of diagram automorphisms
$Q$, given by the image of the vector $\Sigma_{\PP} + \Sigma_{\LL}\xi$ or 
$w_{\PP} +  w_{\LL} \xi \in L \otimes \RR$ where $\xi = ( 1+ i)/\sqrt{2}$
is an eighth root of unity: $\xi^2 = i $. We call this fixed vector the {\it Weyl vector} :
\begin{equation}
\bar{\rho}= (\Sigma_{\PP} + \Sigma_{\LL}\xi)/14
\end{equation}
We note some of the inner products between the special vectors that we need later :
Let $(\rho_1, \dotsb, \rho_{14}) = (x_1, \dotsb, x_7, l_1\xi, \dotsb, l_7\xi)$, 
so that $\sigma$ interchanges $\rho_s$ with $\rho_{7+s}\xi$ and $\bar{\rho}$ is the average
of $\rho_1, \dotsb, \rho_{14}$.
Then we have $\langle \rho_s, \rho_t\rangle $ is equal to $\sqrt{2}$ or 0 according to whether
the two nodes are joined or not joined in the diagram $D$. 
We have, for $ s = 1, \dotsb, 14$,
\begin{equation}
\langle \br, \rho_s \rangle = \nr = 1/(2 + 3\sqrt{2}) 
\end{equation}
From \eqref{wp} and \eqref{wl} we get 
\begin{equation}
\abs{w_{\PP}}^2 = \abs{w_{\LL}}^2 = 2 
\text{  and  } 
\langle \br, w_{\PP} \rangle = \langle \br, w_{\LL} \xi \rangle = 1/\sqrt{2}
\end{equation} 
%
%******************************************************************************************************
%
\begin{heading} The height of a root \end{heading}
We use the Weyl vector $\bar{\rho}$ to define the \textit{height of a root} $r$ as 
\begin{equation*}
\htt(r) = \abs{\langle\bar{\rho}, r\rangle}/\abs{\bar{\rho}}^2 
\end{equation*}
The 14 roots of the diagram $D$ have height equal to 1.
Take the roots of the reflections generating $R(L)$, found in \ref{gen}, given in the coordinate system $\Lambda \oplus H$.
Using the explicit isomorphism found in \ref{calc} we write them in the coordinate system $3E_8 \oplus H$.  
We use the above definition of height and run a ``height reduction algorithm'' (see theorem 5.7 in \cite{tb:el})
on the 81 generators for $R(L)$ found before to see that one can always get to an unit
multiple of an element of $D$. (sometimes one needs  to perturb if the algorithm gets
stuck - at most one perturbation was enough in all cases).
This proves
%
%******************************************************************************************************
%
\begin{theorem}
The order 4 reflections in the roots of $D$ generate the reflection group of $L$.
\end{theorem}
%
%******************************************************************************************************
%
Moreover, as in lemma 3.2 of \cite{tb:el} we can show that the reflections in the roots of the $M_{444}$
diagram generate the reflection group of $L$. We just check that the relation $\mbox{deflate}(y)$ holds
for each octagon $y$ inside $D$. For example, if $y = (d_1,c_1,b_1,a,b_2,c_2,d_2,e_3)$, the relation
$\mbox{deflate}(y)$ follows from the equation
$\phi_{d_1}^{i}\phi_{c_1}^{i}\phi_{b_1}^{i}\phi_{a}^{i}\phi_{b_2}^{i}\phi_{c_2}^{i} (d_2) = -e_3$.
\par
Now we prove the analog of Proposition 6.1 of \cite{tb:el}. The proof is also exactly similar.
\begin{proposition}
The 14 roots of the diagram $D$ are the only roots (up to units) having the  minimum height 1.
All other roots have strictly bigger height.
( In other words, the mirrors of the roots in $D$ are the 14 mirrors closest to the vector $\bar{\rho}$).
\end{proposition}
\begin{proof}
We need the following distance formulae for the metric on the quaternionic hyperbolic space $\HH H^n$
(See \cite{kp:quathyp}).
A positive norm vector $x$ in the vector space determines a point in the hyperbolic space, also denoted by $x$.
A negative norm vector $r$ determines a totally geodesic hyperplane given by $r^{\bot}$. Let
$c(u,v)^2 = \frac{\abs{\langle u,v \rangle}^2}{\abs{u}^2\abs{v}^2}$. Then we have 
\begin{align}
\cosh^2(d(x,x')/2) = c(x,x')^2 , & & \sinh^2(d(x, r^{\bot})/2) = -c(x,r)^2 \label{dp}
\end{align}
Two hyperplanes $ r^{\bot}$ and $r'^{\bot}$ meet in the hyperbolic space if $c(r,r') <1 $, are asymptotic if
$c(r,r') = 1$ and do not meet if $c(r,r') >1$ in which case the distance between the hyperplanes is given by
\begin{equation}
\cosh^2(d(r^{\bot},r'^{\bot})/2) = c(r,r')^2 \label{dhh}
\end{equation} 
Let $r$ be a root of the lattice $L$ with $\htt(r) = \abs{\langle \br, r\rangle} /\nr \leq 1 $. 
We want to prove that $r$ is a unit multiple of one of the 14 roots of $D$.
\par
Let $x$ be a point in $D$. Either $\abs{\langle x, r\rangle} \leq 2$, or
using the triangle inequality $d(r^{\bot}, x^{\bot}) \leq d(r^{\bot}, \br) + d( x^{\bot}, \br) $
along with the distance formulae \eqref{dp} and \eqref{dhh} above we get 
\begin{equation*}
\abs{\langle x, r\rangle} \leq 2\cosh(2\sinh^{-1}(\abs{\br}/\sqrt{2})) \approx 2.32
\end{equation*}
So we must have $\abs{\langle x, r\rangle}^2 $ equal to 0, 2, or 4.
\par
Similarly from $d(r^{\bot}, w_{\PP}) \leq d( r^{\bot}, \br) + d( w_{\PP}, \br)$, \eqref{dp} and \eqref{dhh} we get
\begin{equation*}
\abs{\langle w_{\PP}, r\rangle} \leq 2\sinh(\sinh^{-1}(\abs{\br}/\sqrt{2}) + \cosh^{-1}( 1/2\br) ) \approx 2.26
\end{equation*}
It follows that $\abs{\langle w_{\PP}, r\rangle}^2 $ is equal to 0, 2, or 4.
\par
We can write
\begin{equation}
 r= \sum_{x \in \PP} -x \langle x, r\rangle/2 + w_{\PP} \langle w_{\PP} , r \rangle/2 \label{req}
\end{equation}
taking norm in \eqref{req} we get 
\begin{equation}
  -2 = \sum_{x \in \PP} -\abs{\langle x, r\rangle}^2 /2 + \abs{ \langle w_{\PP} , r \rangle}^2/2 \label{sumip}
\end{equation}
There are only few cases to consider.
Multiplying $r$ by a unit, we may assume that $\langle w_{\PP} , r \rangle$ is either 
0 or $p $ or 2. In the following let $u_1, u_2$ etc denote units in $\cH^*$ and $x_1, x_2$
etc. denote points of $D$.
\par
If $\langle r, w_{\PP}\rangle = 0 $, from \eqref{sumip} we get $\sum \abs{\langle x_s,r\rangle}^2 = 4 $.
Then the unordered tuple $(\langle x_1, r\rangle ,\dotsb, \langle x_{7},r\rangle)$ is equal to
$(2u_1, 0^{6})$ or $((p u_1,p u_2, 0^{5})$.
So either $r$ is an unit multiple of $x_s$ (in which case it has height equal to
one) or $r = (x_1p u_1  + x_2 p u_2  )/ (-2)$. Using
diagram automorphisms (which is 2-transitive on points of $D$) we can assume that
$x_1 = a $ and $x_2 = c_1$ and then check that there is no such root $r$.
\par
If $\langle w_{\PP}, r\rangle = p$, then
$ \sum \abs{\langle x_s,r\rangle}^2 = 6 $; the unordered tuple
$(\langle x_1, r\rangle ,\dotsb, \langle x_{7},r\rangle)$ is equal to 
$(2u_1,  u_2 p, 0^{5})$ or
$(u_1 p,  u_2 p,  u_3 p, 0^4 )$. In the first case
we get $ r = 2x_1 u_1 / (-2) +  x_2 u_2 p / (-2) +  w_{\PP}p/ 2 $.
Taking inner product with $\br$ and using
$\langle \br, w_{\PP} \rangle / \nr = 3 + \sqrt{2}$
we get $\langle \br, r \rangle / \nr = -u_1 - u_2/\bar{p} + (3 + \sqrt{2})/\bar{p} $ 
which clearly has norm greater than one.
\par
In the second case we get 
$ r = \sum_{s = 1}^3  x_s u_s p / (-2) +  w_{\PP} p/ 2 $ which implies
$\langle \br, r \rangle / \nr = (-  u_1 - u_2 - u_3 + 3 + \sqrt{2})/\bar{p} $.
Again this quantity has norm at least one. We now show that the only way it can
be equal to one is if $r$ is a unit multiple of $l_1, \dotsb, l_{7}$.
\par
The only way one can have $\htt(r) = 1$ in the above paragraph
is if $r$ has inner product $p$ with three of the points $x_1, x_2, x_3$
and orthogonal to others. If $x_1,x_2, x_3$ do not all lie on a line
then there is a line $l$ that avoids all these three points. 
Taking inner products with $r$ in the equation
$w_{\PP} = l\bar{p} + \sum_{x \in l} x $ gives
$p = p \langle l,r\rangle $ contradicting $ L'p = L$.
So $x_1, x_2,x_3$ are points on a line $l_1$. It follows that
$r$ and an unit multiple of $l_1$ has the same inner product with each element of
$\PP$ and with $w_{\PP}$. So $r$ is an unit multiple of $l_1$.
\par
If $\langle w_{\PP}, r\rangle = 2 $, and
$\sum \abs{\langle x_s,r\rangle}^2 = 8$, the unordered tuple 
$(\langle x_1, r\rangle ,\dotsb, \langle x_{7},r\rangle)$ is equal to 
$(2u_1, 2u_2, 0^{5})$ or $( 2u_1,  u_2 p,  u_3 p,  0^4)$
or $(u_1 p, \dotsb, u_4 p, 0^3)$. Using similar calculation
as above, we get 
$\langle r, \br\rangle / \nr $ is equal to $(-u_1 -u_2 + 3 + \sqrt{2})$ 
or $(-u_1 -u_2/\bar{p} -u_3/\bar{p} + 3 + \sqrt{2})$ or 
$(-(u_1 + \dotsb + u_4)/\bar{p} + 3 + \sqrt{2})$ respectively.
Again each of these quantities are clearly seen to have norm
strictly bigger than one.
\end{proof}
%
%******************************************************************************************************
%
\begin{heading} Remarks \end{heading}
1. The group generated by the Coxeter diagram $Inc(P^2(\FF_2)$ when
the vertices are made into reflections of order two is the group $O_8^{-}(2):2$
as was found by Simons in \cite{chs:deflate}.
It would be interesting to understand the relation
between these two groups in a more conceptual way.
\par
2. The order of the ``spider element'' $sp = ab_1c_1ab_2c_2ab_3c_3$ in the reflection group $R(L)$
is 40. 
\par
3. Since we can find a $M_{444}$ diagram in $L = \Lambda \oplus H$ with the 6 vectors in the three $E_8$ hands
of the form $(\lambda; 1, *)$ with $\lambda$ in the first shell of Leech lattice, it
is likely that the reflection group of $L$ in-fact equals the whole automorphism group.
It would probably follow after a little more work in the line of Theorem 6.2 of \cite{dja:Leech}
and using the fact that $\TT \subseteq R(L)$. But I have not checked this.
\par 
4. Exactly  as in Remark 3.4 of \cite{tb:el} the $\cH$-lattice $L = 3E_8\oplus H$ can be defined
by the starting with a singular lattice corresponding to the vertices of the diagram $D$ and quotienting out
by the relations  $l(ip) + \sum_{x\in l} x= l'(ip) + \sum_{x\in l'} x $.  
\par
5. The calculations needed for this paper are done using the gp calculator and the
codes (for finding the explicit isomorphism from $3E_8\oplus H $ to $\Lambda\oplus H$ and for the
height reduction algorithm to show that the 14 nodes of $D$ generate $R(L)$)
are contained in the file 
\begin{verbatim}
quat.gp
\end{verbatim}
The programs are available on the my web-site http://www.math.berkeley.edu/\textasciitilde tathagat
%
%*************************************************************************************
%

\end{document}